\documentclass[final, number, 3p, 12pt]{elsarticle}
\biboptions{sort&compress}

\usepackage{epsfig}

\usepackage{amssymb}
\usepackage{amsmath}
\usepackage{hyperref}
\usepackage{float}

\usepackage{svg}

\newcommand{\fig}[1]{\figurename~\ref{#1}}
\newcommand{\tbl}[1]{\tablename~\ref{#1}}

\begin{document}

\begin{frontmatter}



\title{Enhanced Diffuse Interface Method for Multiphase Flow Simulations Across All Mach Numbers}


\author[label1]{Ghanshyam Bharate}
\author[label1]{J.C. Mandal}
\affiliation[label1]{organization={Department of Aerospace Engineering, Indian Institute of Technology Bombay},  city={Mumbai}, postcode={400076}, country={India}}

\begin{abstract}
 This paper enhances the Diffuse Interface Method (DIM) for simulating compressible multiphase flows across all Mach numbers by addressing the accuracy challenges posed at low Mach regimes. A correction to the Riemann solver is introduced, designed to mitigate excessive numerical diffusion while maintaining simplicity and efficiency. The validity of this correction is established through rigorous asymptotic analysis of the governing equations and their discrete counterparts. The proposed correction is implemented within a six-equation model framework with instantaneous relaxation using an HLLC-type solver. Numerical test cases demonstrate significant improvements in accuracy, confirming the effectiveness of the approach in capturing multiphase flow dynamics across a wide range of Mach numbers. 
\end{abstract}

\begin{keyword}
Diffused Interface Method (DIM); Six-equation Model; Multiphase Flows; All Mach Number; HLLC Solver



\end{keyword}

\end{frontmatter}

\section{Introduction}

Accurate numerical simulation of multiphase flows involving liquid and gas is essential for a wide range of industrial and scientific applications. While flows such as cavitation, boiling, and depressurization in transport systems are generally classified as low-Mach-number regimes, they can exhibit localized compressible regions. Capturing these regions accurately requires a compressible model and solver. Furthermore, the significant disparity in the speed of sound between liquid and gas phases underscores the necessity of an all-Mach compressible framework to ensure robust and precise computations.

In recent years, diffuse interface methods (DIMs) have emerged as highly effective tools for simulating compressible multiphase flows, including dispersive multiphase mixtures and interfaces between pure fluids. Treating fluid-fluid interfaces as contact discontinuities ensures precise wave transmission. Furthermore, their unique ability to dynamically generate material interfaces not initially present has contributed significantly to their growing popularity. 

There are various models of DIMs existing in literature \cite{saurel1999, kapila2001, Saurel2009SimpleAE} . The most general form is the seven-equation model. The seven-equation is unconditionally hyperbolic and fully non-equilibrium. This model is first introduced by \citet{baer1986} for solid combustible granular flows and then it is modified by \citet{saurel1999} for the computation of general multiphase flow problems. They \cite{saurel1999} have also included relaxation terms in the seven equation model which takes drag and pressure relaxation effects into account. With the help of instantaneous relaxation procedure this fully non-equilibrium model can be used for the computation of equilibrium flows. 

 Despite the numerous features of the seven-equation model, it leads to a large system with a substantial number of waves. This has prompted researchers to explore alternative computationally inexpensive models, such as five-equation model \cite{murrone2005five} and six-equation model \citep{Saurel2009SimpleAE}. These simplified models are derived in the zero relaxation time limit from the seven-equation model. They were proposed by \citet{kapila2001} for granular energetic flows. Although the five equation model \cite{murrone2005five} is computationally cheaper among all the discussed models, it is also associated with several numerical challenges \citep{Saurel2009SimpleAE, petitpas2007relaxation}. These challenges include maintaining volume fraction positivity and dealing with the non-monotonic speed of sound. Because of these difficulties, several researchers \cite{Saurel2009SimpleAE, zein2010, pelanti2014mixture, yu2023numerical} in the past have adopted the approach where a two-pressure single-velocity six-equation model is solved with instantaneous relaxation procedure for the computation of the mechanical equilibrium (single-velocity single-pressure) flows. In this work, the same six-equation model is used for the computation of multiphase flows.

 The shock capturing ability of the Riemann solvers in supersonic and transonic flows makes it an ideal choice for computing numerical fluxes at cell interfaces. However, the inherent numerical viscosity, essential for maintaining stability in the presence of strong discontinuities, becomes problematic in low Mach limit $(M \rightarrow 0)$. In such cases, the numerical viscosity causes excessive diffusion, resulting in unphysical outcomes and incorrect pressure scaling. The detailed explanation of the unphysical behaviour of the Riemann solvers in low Mach range for single phase flow can be found in the literature \cite{guillard1999behaviour, rieper2011low, dellacherie2010analysis}.

In the past many researchers \cite{guillard1999behaviour, guillard2004behavior, rieper2011low, osswald2016l2roe, luo2005extension, pelanti2018wave, HLLC_IJNMF, gogoi2023low, gogoi2024simple, gogoi2025enhanced} have presented different methods for addressing the unphysical low Mach number behaviour of the Riemann solvers for single phase flow. These methods can be broadly classified into two categories. 

The first category comprises preconditioning methods \cite{guillard1999behaviour, guillard2004behavior, luo2005extension, pelanti2018wave}, where a preconditioned Riemann problem is solved using an existing Riemann solver. In this approach, the modified wave speeds computed through complex algebraic expressions are used for the computation of fluxes. While, the preconditioning methods help in reducing the excessive numerical viscosity and rectify the incorrect scaling of pressure, they have notable drawbacks. The primary issue is the global cut-off Mach number problem \cite{li2008all}. Additionally in explicit scheme, these methods impose a severe time step restriction, requiring time step size proportional to the square of the Mach number $(\Delta t  \propto M^2 )$ \cite{birken2005stability}. 

The second category of methods consists of correction methods developed by researchers \cite{thornber2008improved, rieper2011low, osswald2016l2roe,  HLLC_IJNMF, gogoi2023low, gogoi2024simple, gogoi2025enhanced}. These methods address low Mach number issues in Riemann solvers by scaling the velocity jumps, providing a simple and straightforward implementation. Unlike preconditioning methods, they do not suffer from global cut-off problem, as they rely on the local Mach number to scale the velocity jump. Moreover, they impose no severe time step restriction.

Attempts \cite{murrone2008, liquidgas2013, pelanti2017low} have also been made to extend Riemann solvers to low Mach number multiphase flows. \citet{murrone2008} have adapted the single phase preconditioning approach \cite{guillard2004behavior} to the five equation multiphase model \cite{murrone2005five}. Similarly, other researchers \cite{liquidgas2013, pelanti2017low} applied this approach to six-equation models \cite{Saurel2009SimpleAE, pelanti2014mixture}. However, all of these efforts fall within the preconditioning category.

In this work, we focus on a different class of methods based on low Mach correction that have not yet been implemented in multiphase Riemann solvers. Our objective is to develop a numerical scheme capable of computing multiphase flows across all Mach numbers using the six-equation model \cite{Saurel2009SimpleAE}. For the flux computation, we utilize an HLLC-type Riemann solver, and implement the correction method suggested by \citet{thornber2008improved} to address the low Mach number challenges.

The paper is structured into ten sections. Multiphase model is described in Section \ref{sec:model}. The six-equation model is solved in two steps: evolution and relaxation~\cite{Saurel2009SimpleAE, zein2010, pelanti2014mixture}. Both steps are discussed in Sections \ref{sec:evol_step} and \ref{sec:relax_step}. Higher order formulation used in this work is presented in the Section \ref{sec:second_order}. The behaviour of the continuous model in the low Mach limit $(M \rightarrow 0)$, analysed using asymptotic expansion, is discussed in Section \ref{sec:asympt_contin}. The correction proposed to address the low Mach problem in the existing numerical method is presented Section \ref{sec:low_mach}. The effectiveness of the proposed correction is demonstrated in Section \ref{sec:asympt_discrete} by comparing asymptotic expansion of the discretised equations with and without low Mach correction. The results of the several test cases are presented and discussed in the Section \ref{sec:results}. Section \ref{sec:conclusion} presents the conclusion of the paper.

\section{Multiphase model}\label{sec:model}
The six-equation multiphase model \cite{Saurel2009SimpleAE} can be written as

\begin{equation}\label{eq:model}
 \begin{gathered}
      \frac{\partial\mathbf{U}}{\partial t} +  \nabla \cdot {H} (\mathbf{U})  + \mathbf{\sigma}\left(\mathbf{U}\right)  =   \mathbf{S} (\mathbf{U}) \\ 
      {H} (\mathbf{U}) = (\mathbf{F}, \mathbf{G}) \\
      \mathbf{U} = \left[\begin{array}{c} \alpha_1 \\ \alpha_1 \rho_1 \\ \alpha_2 \rho_2 \\  \rho u \\ \rho v \\ \rho E \\\alpha_1 \rho_1 e_1 \\ \alpha_2 \rho_2 e_2 \end{array}\right], \quad
     \mathbf{F} = \left[\begin{array}{c} 0 \\ \alpha_1 \rho_1 u \\  \alpha_2 \rho_2 u  \\  \rho {u^2} +  p  \\ \rho u v \\ (\rho E + p) u \\ \alpha_1 \rho_1 e_1 u \\ \alpha_2 \rho_2 e_2 u \end{array}\right], \quad
     \mathbf{G} = \left[\begin{array}{c} 0 \\ \alpha_1 \rho_1 v \\  \alpha_2 \rho_2 v  \\  \rho u v  \\  \rho v^2 + p \\ (\rho E + p) v \\ \alpha_1 \rho_1 e_1 v \\ \alpha_2 \rho_2 e_2 v \end{array}\right]\\
    \mathbf{\sigma}\left(\mathbf{U}\right) = \left[\begin{array}{c} \mathbf{u} \cdot \nabla \alpha_1 \\ 0 \\ 0 \\ 0 \\ 0 \\ 0 \\  \alpha_1 p_1 \nabla  \cdot \mathbf{u} \\ \alpha_2 p_2 \nabla \cdot \mathbf{u} \\ \end{array}\right], \quad
    \mathbf{S}(\mathbf{U}) = \left[\begin{array}{c}  \mu \left(p_1 - p_2 \right)  \\ 0 \\ 0 \\0\\ 0 \\ 0\\ -\mu p_I \left(p_1 - p_2 \right) \\ \mu p_I \left(p_1 - p_2 \right)\end{array}\right]
      \end{gathered}
\end{equation}
where $\mathbf{U}$ is the set of conservative variables and ${H}$ is a conservative flux tensor. $\mathbf{\sigma}\left(\mathbf{U}\right)$ and $\mathbf{S}$ contains non-conservative and pressure relaxation terms.

Here $\alpha_j, \rho_j, p_j, e_j$ represent the volume fraction, density, pressure, and specific internal energy of the phase $j$, respectively. The velocity vector is represented by $\mathbf{u} = (u, v)$, while $\rho$ and $p$ stands for the mixture density and mixture pressure. The expressions for these mixture quantities is given by

\begin{equation}
 \quad p = \sum_j \alpha_j p_j, \quad \rho = \sum_j \alpha_j \rho_j
\end{equation} 
Interface pressure $(p_I)$ appearing in the relaxation terms is defined as 

\begin{equation}
p_I = \frac{z_1 p_2 + z_2 p_1}{z_1 + z_2}
\end{equation}
here, $z_j = \rho_j a^2_j$ represents acoustic impedance of phase $j$. The expression for mixture total energy $(\rho E)$ can be written as

\begin{equation}
\rho E = \rho e +   \frac{\rho \mathbf{u} \cdot \mathbf{u}}{2}
\end{equation}
The relation between phasic internal energy $(e_j)$ and mixture internal energy $(e)$ is given by following expression.
\begin{equation}
\rho e  = \sum_j \alpha_j \rho_j e_j
\end{equation}

Internal energy $(e_j)$ and speed of sound $(a_j)$ are determined using following stiffened gas equation of state (SGEOS) relations
\begin{equation}\label{eq:EOS}
e_j  = \frac{p_j + \gamma_j \pi_j}{\rho_j \left(\gamma_j -1 \right)}, \quad a_j = \sqrt{\frac{\gamma_j \left(p_j + \pi_j \right)}{\rho_j}}
\end{equation}

The volume fraction $(\alpha_j)$ for the multiphase system can be defined as fraction of total volume occupied by phase $j$ and it always comply with following saturation condition

\begin{equation}
\sum_j \alpha_j = 1
\end{equation}  
In order to avoid numerical problems like infinite density and pressure, the $\alpha_j$ should be always non-zero quantity, hence in absence of particular fluid it can be set to a very low value $\alpha_j = \epsilon$, where $\epsilon$ represent very small number.  

The multiphase model described in \eqref{eq:model} is an overdetermined system, that includes an extra mixture total energy equation. The extra equation was introduced by~\citet{Saurel2009SimpleAE} to ensure conservation of mixture total energy. The homogeneous part of six-equation model is a hyperbolic system and its one-dimensional version has six real eigenvalues, given as   
\begin{equation}
\lambda_1 = u - a, \quad \lambda_{2,..,5} = u, \quad \lambda_6 = u + a.
\end{equation}
Here, $a$ represents mixture sound speed, which can be computed as  
\begin{equation}\label{eq:mixture_sound}
a = \sqrt{\frac{1}{\rho}  \sum_j \alpha_j \rho_j a^2_j}
\end{equation} 

The multiphase system \eqref{eq:model} with relaxation source terms is solved using operator splitting approach. This method includes a two step procedure to obtain the solution for each time level. The first step is an evolution step in which equations \eqref{eq:model} are solved without any relaxation terms using a hyperbolic operator $\mathcal{L}_H$. The second step is a relaxation steps, which takes account for only relaxation effects and this is done by solving system of ODEs: $\mathbf{U}_t = \mathbf{S(U)}$ using relaxation operator $\mathcal{L}_R$. In this study, our focus is on the computation of equilibrium flows, hence the relaxation step is performed in the limit  $\mu \rightarrow \infty$ and by doing so we are instantaneously relaxing pressures of both phases to the common value. The solution at $n+1$ time level in terms of operators can be written as

\begin{equation}
\begin{gathered}
\mathbf{U}^{n+1} = {\mathcal{L}_R\mathcal{L}_H} \left( \mathbf{U}^{n} \right) \\
\mathcal{L}_H: \frac{\partial\mathbf{U}}{\partial t}+ \nabla \cdot \mathcal{F}  + \mathbf{\sigma}\left(\mathbf{U}\right)  = 0 \\
\mathcal{L}_R: \frac{\partial\mathbf{U}}{\partial t} = \mathbf{S}
\end{gathered}
\end{equation}

\begin{figure}[H]
           \centering
             \includegraphics[ scale = 1]{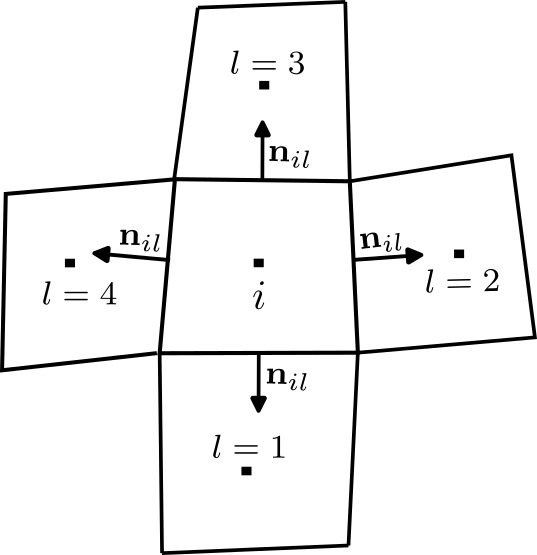}
           \caption{Stencil of quadrilateral cells with notation}
           \label{fig:fvm}
\end{figure}

\section{Evolution step}\label{sec:evol_step}
In this section we look at the solution of multiphase model \eqref{eq:model} without any relaxation terms. The multiphase system \eqref{eq:model} contain non-conservative terms in volume fraction and phasic internal energy equations. These terms required special attention for discretization. As per the discretization method, the homogeneous part of multiphase system \eqref{eq:model} can be divided into two groups. One group contains conservative part of the system which includes mass, momentum and total energy equation. Another group has remaining equations, which have non-conservative terms. First, we look at the discretization of conservative group of equations, these system of equations can be written as  

\begin{equation}\label{eq:model_conservative}
    \begin{gathered}
     \frac{\partial\mathbf{Q}}{\partial t}+ \nabla \cdot {H} (\mathbf{Q})    =   0 \\
     {H} (\mathbf{Q})   = (\mathbf{F}, \mathbf{G})
    \end{gathered}
    \end{equation}
    where 
    \begin{equation}
    \begin{gathered}
    \mathbf{Q} = \left[\begin{array}{c} \alpha_1 \rho_1 \\ \alpha_2 \rho_2 \\  \rho u \\ \rho v \\  \rho E \end{array}\right], \quad \mathbf{F}  = \left[\begin{array}{c}  \alpha_1 \rho_1 u \\  \alpha_2 \rho_2 u  \\  \rho {u^2} +  p  \\ \rho u v \\  (\rho E + p) u \end{array}\right], \quad
     \mathbf{G} = \left[\begin{array}{c}  \alpha_1 \rho_1 v \\ \alpha_2 \rho_2 v\\  \rho u v \\ \rho {v^2} + p  \\ (\rho E + p) v  \end{array}\right]
     \end{gathered}
    \end{equation}

The finite volume discretization of above system for quadrilateral cell is given by 

\begin{equation}\label{eq:FVM}
    \begin{gathered}
        \frac{\mathrm{d} \mathbf{Q_i}}{\mathrm{d} t}  + \frac{1}{\Omega_i} \sum_{l = 1}^{4}  H_{il}  \cdot  \mathbf{n}_{il} ~{\Delta s}_{il} = 0
    \end{gathered}
\end{equation}

Here, ${\Delta s}_{il}$ represents the length and $\mathbf{n}_{il}$ represents the normal face vector of cell boundary between cell $i$ and its neighbouring cell $l$. The stencil of neighbouring cells with notation is shown in \fig{fig:fvm}. Using rotation invariance property the normal flux $(H_{il}  \cdot  \mathbf{n}_{il})$ can be written as

\begin{equation}\label{eq:rotational1}
    H_{il}  \cdot  \mathbf{n}_{il} = \mathrm{T}^{-1}_{il} \mathbf{F} \left( \mathbf{\hat{Q}}_{L}, \mathbf{\hat{Q}}_{R} \right), 
  \end{equation}

where $\mathrm{T}^{-1}_{il}$ is the inverse rotational matrix and $ \mathbf{F} \left( \mathbf{\hat{Q}}_{L}, \mathbf{\hat{Q}}_{R} \right)$ is the normal conservative flux vector in locally rotated coordinate, which can be computed by the Riemann solver.  $\mathbf{\hat{Q}}_{L}$  and $\mathbf{\hat{Q}}_{R}$ are set of normal conservative variables at left and right side of the interface in the locally rotated coordinate system. The normal variables can be computed as

\begin{equation}\label{eq:rotational2}
    \begin{gathered}
      \mathbf{\hat{Q}}_{L/R} =  \mathrm{T}_{il} \mathbf{Q}_{L/R} = \left[\begin{array}{c}
    \alpha_1 \rho_1 \\ \alpha_2 \rho_2 \\  \rho {u}_n \\ \rho u_t \\ \rho E
    \end{array}\right]_{L/R} , \quad \mathrm{T}_{il} = \left[\begin{array}{ccccc}
    1 & 0 & 0 & 0 & 0\\
    0 & 1 & 0 & 0 & 0\\
    0 & 0 & n_{x} & n_{y} & 0\\
    0 & 0 &-n_{y} & n_{x} & 0\\
    0 & 0 & 0 & 0 & 1 
    \end{array}\right] 
    \end{gathered}
\end{equation}
Here, ${u}_n$ and ${u}_t$ are normal and tangential velocities at cell interface, and $(n_{x}, n_{y})$ are components of unit normal face vector $\mathbf{n}_{il}$. 

For volume fraction and phasic internal energy equations, semi-discrete formulations can be written as

\begin{equation}\label{eq:non_conserv}
    \begin{gathered}
 \frac{d{(\alpha_1)}_i}{dt} + \frac{1}{\Omega_i} \sum_{l = 1}^{4} \left[ {\left(\alpha_1 u_n \right)}_{il} - {(\alpha_1)}_i ({u}_{n})_{il} \right] ~ \Delta s_{il} = 0 \\
    \frac{d (\alpha_j \rho_j e_j)_i}{dt} + \frac{1}{\Omega_i} \sum_{l = 1}^{4} \left[{\left( \alpha_j \rho_j e_j u_n \right)}_{il} +    {(\alpha_j p_j)}_i ~ ({u}_{n})_{il}\right] \Delta s_{il}  = 0 , \quad j = 1, 2
\end{gathered}
\end{equation}

In the formulation presented in equation \eqref{eq:non_conserv}, the terms with subscript $i$ are computed using cell-averaged variables. In contrast, the quantities with subscript $il$ are obtained at the cell interface using an HLLC-type Riemann solver \cite{Saurel2009SimpleAE}, as explained below.

\subsection{HLLC type solver}

The expression to determined normal flux vector, $ \mathbf{F} \left( \mathbf{\hat{Q}}_{L}, \mathbf{\hat{Q}}_{R} \right)$ for the HLLC solver can be written as

\begin{equation}\label{eq:HLLC}
\begin{gathered}
\mathbf{F} \left( \mathbf{\hat{Q}}_{L}, \mathbf{\hat{Q}}_{R} \right) =  \left\lbrace \begin{array}{lll} 
\mathbf{F}\left( \mathbf{\hat{Q}}_{L} \right), & if & 0 \leq S_L \\
\mathbf{F}\left( \mathbf{\hat{Q}}_{L} \right) + S_L \left(\mathbf{\hat{Q}}^{*}_{L} - \mathbf{\hat{Q}}_{L} \right) , & if & S_L \leq 0 \leq S^{*} \\
\mathbf{F}\left( \mathbf{\hat{Q}}_{R} \right) + S_R \left(\mathbf{\hat{Q}}^{*}_{R} - \mathbf{\hat{Q}}_{R} \right) , & if & S^{*} \leq 0 \leq S_R \\
\mathbf{F}\left( \mathbf{\hat{Q}}_{R} \right) , & if & S_R \leq 0
\end{array}  \right.
\end{gathered}
\end{equation}

Where, $S_L$, $S_R$, and $S^{*}$ are speeds of left running wave, right running wave and contact discontinuity. The set of intermediate conservative variables $\mathbf{\hat{Q}}^{*}_{K}$ near $S_K$ wave can be written as

\begin{equation}\label{eq:Q_k}
\mathbf{Q}^{*}_K =  \left(\frac{S_K - {u}_{nK}}{S_K - {S}^{*}} \right)  \left[ \begin{array}{c} { \alpha_{1K} \rho_{1K} } \\ { \alpha_{2K} \rho_{2K} }  \\ \rho_{K} {S}^{*} \\ \rho_{K} u_{tK} \\  \rho_{K}\left( E_{K} + \left(S^{*} - {u}_{nK} \right) \left( S^{*} + \dfrac{p_{K}}{\rho_{K}\left( S_K - {u}_{nK} \right)} \right) \right)
\end{array} \right], \quad K = L~\text{or}~R
\end{equation}  
where speed of contact wave $S^{*}$ is given by following expression: 

\begin{equation*}
S^{*} = \frac{  p_L - p_R +   {\rho}_{R} {u}_{nR} (S_R - {u}_{nR}) - \rho_L {u}_{nL} (S_L - {u}_{nL}) }{ {\rho}_R (S_R - {u}_{nR}) - {\rho}_L (S_L - {u}_{nL})}
\end{equation*}

The above expression is identical to the expression given by \citet{riemann1997} for single phase HLLC solver, however here mixture quantities $(p_{K}, \rho_{K})$  are used instead single phase quantities. The wave speeds, $S_L$ and $S_R$ can be estimated using following expressions \cite{davis1988simplified}:

\begin{equation}
S_R = \max({u}_{nL} + a_L, {u}_{nR} + a_R), \quad S_L = \min({u}_{nL} -a_{L}, {u}_{nR} - a_R) 
\end{equation}

For semi-discrete equations \eqref{eq:non_conserv}, we also require phasic quantities across left and right waves. Such as, volume fraction $(\alpha^{*}_{jK})$, internal energy $(e^{*}_{jK})$ and density $(\rho^{*}_{jK})$. As the volume fraction is constant across wave $S_{K}, K = L~\text{or}~R$,  $\alpha^{*}_{jK}$ can be simply taken as \cite{Saurel2009SimpleAE}
\begin{equation}
\alpha^{*}_{jK} = \alpha_{jK}.
\end{equation}
From \eqref{eq:Q_k} we can write expression of phasic density $(\rho^{*}_{jK})$ as 
\begin{equation}
\rho^{*}_{jK} = \rho_{jK} \frac{S_K - u_{nK}}{S_K - S^{*}}.
\end{equation}
The phasic internal energy $(e^{*}_{jK})$ is determined from SGEOS relation~\eqref{eq:EOS}, which requires phasic pressure $( p^{*}_{jK})$ across wave $S_K$ and it can estimated using following expression \cite{Saurel2009SimpleAE} 

\[ p^{*}_{jK} = (p^{*}_{jK} + \pi_j) \frac{(\gamma_j -1)\rho_{jK} - (\gamma_j + 1)\rho^{*}_{jK} }{(\gamma_j -1)\rho^{*}_{jK} - (\gamma_j + 1)\rho_{jK}} - \pi_j \]

\section{Relaxation step}\label{sec:relax_step}
In the relaxation step, the non-equilibrium solution from the evolution step is brought to a state of mechanical equilibrium using an instantaneous pressure relaxation procedure. In this step, the following system of ODEs is solved in the limit  $\mu \rightarrow \infty$. 
\begin{equation}\label{eq:relax_ODE}
\begin{gathered}
 \frac{\partial\mathbf{U}}{\partial t} = \mathbf{S} \\
   \mathbf{U} = \left[\begin{array}{c} \alpha_1 \\ \alpha_1 \rho_1 \\ \alpha_2 \rho_2 \\  \rho u \\ \rho v \\ \rho E \\ \alpha_1 \rho_1 e_1 \\ \alpha_2 \rho_2 e_2 \end{array}\right], \quad
\mathbf{S} = \left[\begin{array}{c}  \mu \left(p_1 - p_2 \right) \\ 0 \\ 0 \\ 0 \\ 0 \\0 \\-\mu p_{I} \left(p_1 - p_2 \right) \\ \mu p_I \left(p_1 - p_2 \right)\end{array}\right]
 \end{gathered}
\end{equation}
As we can see from the system of equations \eqref{eq:relax_ODE}, the following quantities remain unchanged:   

\begin{equation}\label{eq:relax_constant}
{(\alpha_j \rho_j)}^{*} = {(\alpha_j \rho_j)}^{o}, \quad j = 1,2, \quad  \rho^{*} = \rho^{o} \quad \text{and} \quad \mathbf{u}^{*} = \mathbf{u}^{o}.
\end{equation}
Here, the superscript $"o"$ and $"*"$ is used for the variables before and after relaxation step. Using relations \eqref{eq:relax_constant} and volume fraction equation \eqref{eq:relax_ODE} in phasic internal energy equations \eqref{eq:relax_ODE}, we obtain following system of ODEs.  

\begin{equation}\label{eq:relax_ODE2}
\begin{gathered}
\frac{\partial {\left( \alpha_j \rho_j e_j \right)}}{\partial t} = - p_I \frac{\partial {\alpha_j }}{\partial t}, \quad j = 1, 2
\end{gathered}
\end{equation}
The numerical discretization of above ODEs \eqref{eq:relax_ODE2} can be written as

\begin{equation}\label{eq:relax_eq2}
\begin{gathered}
 {\left( \alpha_j \rho_j e_j \right)}^{*} - {\left( \alpha_j \rho_j e_j \right)}^{o} = - \bar{p}_I {\left( \alpha^{*}_j  - \alpha^{o}_j \right) } , \quad j = 1,2
\end{gathered}
\end{equation}
Here, $\bar{p}_I$ represents the numerical approximation of $\frac{1}{\left( \alpha^{*}_j  - \alpha^{o}_j \right)} \int p_I \partial {\alpha_j } $ and it can be taken as weighted average of final and initial values. 

\begin{equation}\label{eq:aver_PI}
 \bar{p}_I = (1 - \omega) p^{o}_I + \omega p^{*}_I, \quad \omega \in [0, 1]  
\end{equation}
There are three possible choices: $\omega = 0$, $\omega = 1,$ and $\omega = \frac{1}{2}$. However, no significant difference is observed in the results obtained by these choices \cite{Saurel2009SimpleAE}. Here we have chosen $\omega = \frac{1}{2}$.

After substituting $\rho_j e_j$ with $p_j$ using the SGEOS relation \eqref{eq:EOS} and applying the mechanical equilibrium condition, $p^{*}_1 = p^{*}_2 = p^{*}_I = p^{*}$, the number of variables in final equations \eqref{eq:relax_eq2} are reduced to two, namely $p^{*}$ and $\alpha^{*}_1$. By solving these two equations  with two unknowns we can easily obtain solution after the relaxation step.

\subsection{Reinitialization step}

The phasic pressure $p^{o}_j$ and the phasic internal energy $e^{o}_j$ prior to the relaxation step are determined by solving non-conservative internal energy equations. Since these equations do not ensure
the conservation of mixture total energy, variables after relaxation step may also violate the conservation principle. To ensure conservation of mixture total energy, Saurel
et al. [22] introduced the idea of Reinitialization of phasic pressure and phasic energy. The reinitialized pressure, $p^{**}$ is computed using the following expression.

\begin{equation}
\begin{gathered}
 p^{**} = \dfrac{(\rho e)^{o} - \sum\limits_{j} \alpha^{*}_j \left( \dfrac{\gamma_j \pi_j}{\gamma_j -1} \right)}{ \sum\limits_{j} \dfrac{\alpha^{*}_j}{(\gamma_j -1)}}. 
 \end{gathered}
\end{equation} 

Here, $(\rho e)^{o}$ is computed from mixture total energy equation, and $\alpha^{*}_j$ is taken from the solution obtained after relaxation step. After getting reinitialized pressure $p^{**}$, internal energies are also reset according to the $p^{**}$.


\section{Second order formulation}\label{sec:second_order}

To achieve second-order spatial accuracy, the set of primitive variables at the left and right state, $\mathbf{W}_{L/R} = \left[\alpha_1, \rho_1, \rho_2, u, v, p_1, p_2\right]_{L/R}$, are reconstructed using the cell average values $(\mathbf{W}_{i}, \mathbf{W}_{l})$. The reconstruction formulae using a truncated Taylor series expansion can be written as

\begin{equation}\label{eq:reconstruction}
  \begin{gathered}
  W_L = W_{i} + \left[  \left(\frac{\partial W}{\partial x}\right)_{i} (x_{il} - x_i) + \left(\frac{\partial W}{\partial y}\right)_{i} (y_{il} - y_i)\right] \\
  W_R = W_{l} + \left[  \left(\frac{\partial W}{\partial x}\right)_{l} (x_{il} - x_{l}) + \left(\frac{\partial W}{\partial y}\right)_{l} (y_{il} - y_l) \right]
  \end{gathered}
\end{equation}

Here, $(x_i, y_i)$, $(x_{l}, y_{l})$ and $(x_{il}, y_{l})$ are the coordinates of centroid of cell $i$, neighbouring cell $l$, and the centre of cell face $il$. Using the same formulae in \eqref{eq:reconstruction}, following linear system of equations is formulated for computing derivatives, ${\left(\frac{\partial W}{\partial x}\right)}_i$, ${\left(\frac{\partial W}{\partial y}\right)}_i$.

\begin{equation}
\begin{gathered}
 \underbrace{ \left[\begin{array}{ll}
    (x_1 -x_i) &  (y_1 - y_i) \\
    (x_2 -x_i) & (y_2 - y_i) \\
    (x_3 -x_i) & (y_3 - y_i) \\
    (x_4 -x_i) & (y_4 - y_i) 
    \end{array}\right]}_{\mathrm{S}}
     \underbrace{  \left[\begin{array}{l}
        {\left(\dfrac{\partial W}{\partial x}\right)}_i \\
        {\left(\dfrac{\partial W}{\partial y}\right)}_i \\
        \end{array}\right]}_{\mathbf{dW}} = \underbrace{ \left[\begin{array}{l}
   W_1 - W_i\\
   W_2 - W_i\\
   W_3 - W_i\\
   W_4 - W_i
  \end{array}\right]}_{\mathbf{\Delta W}}
\end{gathered}
\end{equation}

This overdetermined system is solved using with SWDLS \cite{SDWLS, MANDAL2015669} method. This method utilizes the weighted least-square approximation, which can be written as 

\begin{equation}
\mathbf{dW} = {\left(\mathrm{S}^{T} \mathrm{w} \mathrm{S}\right)}^{-1} \mathrm{S}^{T} \mathrm{w} \mathbf{\Delta W}
\end{equation}

The matrix, $\mathrm{w} = diag(w_1,.., w_7)$ contain, solution dependent weights, $w_l = \frac{1}{\Delta W^2_l + \epsilon}$. Here, $\epsilon$ is a very small number. For second order accuracy in time, SSPRK method is used. The steps for implementing SSPRK scheme are
\begin{equation}
\begin{split}
\mathbf{U}^{(1)} &= \mathcal{L}^{\Delta t}_R \mathcal{L}^{\Delta t}_H \left( \mathbf{U}^{n} \right)  \\ 
\mathbf{U}^{(2)} &= \mathcal{L}^{\Delta t}_R \mathcal{L}^{\Delta t}_H\left(\mathbf{U}^{(1)}\right) \\
\mathbf{U}^{(3)} &= \frac{1}{2}\left( \mathbf{U}^n + \mathbf{U}^{(2)} \right)  \\
\mathbf{U}^{n+1} &= \mathcal{L}^{\Delta t}_R \left(\mathbf{U}^{(3)}\right) 
\end{split}
\end{equation}

\section{Asymptotic analysis of continuous model}\label{sec:asympt_contin}
In order to understand the low Mach number problem associated with the Riemann solver, researchers~\citep{rieper2011low, osswald2016l2roe, murrone2008, pelanti2017low} in the past have done asymptotic analysis of continuous equations as well as semi-discrete equations. The asymptotic expansion of non-dimensionalized equations reveals the behaviour of flow variables in the low Mach number limit. In the case of diffuse interface methods, \citet{murrone2008} were the first to report the asymptotic analysis of five-equation model. Later \citet{liquidgas2013} and then \citet{pelanti2017low} have done the  asymptotic analysis for six-equation model.  Since we are interested in solving six-equation model, we will briefly discuss the asymptotic analysis of later model. After non-dimensionalization, the homogeneous part of equations \eqref{eq:model} can be written as
\begin{equation}\label{eq:non_dimensional_equation}
 \begin{split}
  & \frac{\partial \tilde{\alpha}_1}{\partial \tilde{t}} + \tilde{\mathbf{u}} \cdot \nabla \tilde{\alpha}_1  = 0  \\
  &  \frac{\partial \tilde{\alpha}_j  \tilde{\rho}_j }{\partial \tilde{t}} + \nabla \cdot \left(\tilde{\alpha}_j \tilde{\rho }_j \tilde{\mathbf{u}} \right) = 0 \\
  &  \frac{\partial \tilde{\rho } \tilde{\mathbf{u}} }{\partial \tilde{t}} + \nabla \cdot \left( \tilde{\rho} \tilde{\mathbf{u}} \otimes \tilde{\mathbf{u}} \right) + \frac{1}{M^2} \nabla \tilde{p} = 0 \\
  &  \frac{\partial \tilde{\rho } \tilde{E }}{\partial \tilde{t}} + \nabla \cdot \left( \left( \tilde{\rho } \tilde{E } + \tilde{p} \right) \tilde{\mathbf{u}} \right) = 0  \\
      &   \frac{\partial  \tilde{\alpha}_j  \tilde{\rho}_j  \tilde{e}_j }{\partial \tilde{t}} + \nabla \cdot \left(\tilde{\alpha}_j  \tilde{\rho}_j  \tilde{e}_j  \tilde{\mathbf{u}} \right) + \tilde{\alpha}_j  \tilde{p}_j  \nabla \cdot \tilde{\mathbf{u}}  = 0, \quad j = 1, 2
\end{split}
\end{equation}
The non-dimensionalization of the variables are performed in the following manner:
\begin{equation}\label{eq:reference}
\begin{split}
& \tilde{\alpha}_j = \alpha_j, \quad  \tilde{\rho}_j = \frac{\rho_j}{[\rho]}, \quad \tilde{p}_j = \frac{p_j}{[\rho][a]^2},\quad j = 1, 2 \\
& \tilde{\mathbf{u}} = \frac{\mathbf{u}}{[u]}, \quad  \tilde{\mathbf{x}} = \frac{\mathbf{x}}{[x]}, \quad \tilde{t}  = t \frac{[u]}{[x]}
  \end{split}
\end{equation}
For any variable ${\phi}$, $\tilde{\phi}$ represents the variable after non-dimensionalization and $[\phi]$ represents reference quantity used for non-diemnsionalization. It can be seen the non-dimensional momentum equation \eqref{eq:non_dimensional_equation}, there is an extra factor $\frac{1}{M^2}$ appearing before mixture pressure gradient term. Here, $M$ is reference Mach number, which can be defined as: $M = \frac{[u]}{[a]}$. For non-equilibrium six-equation model, $M$ is based on the reference mixture sound speed $([a])$. Any variable $\tilde{\phi}$ can be expanded in terms of reference Mach number as   

\begin{equation}\label{eq:asym_expansion}
 \tilde{\phi} = \tilde{\phi}^{(0)} + \tilde{\phi}^{(1)} M + \tilde{\phi}^{(2)} M^2  
\end{equation}

After substituting above expansion for the variables appearing in \eqref{eq:non_dimensional_equation} we get separate equations according to the order of $M$. Since all quantities are non-dimensionalized, we have removed the superscript  $\tilde{(\hspace{2pt})}$ for brevity in the following asymptotic results.

\begin{enumerate}
\item Order of $M^{-2}$

\begin{equation}\label{eq:count_M_2}
  \nabla {p}^{(0)} = 0
\end{equation}

\item Order of $M^{-1}$ terms

\begin{equation}\label{eq:count_M_1}
   \nabla {p}^{(1)} = 0
\end{equation}

\item Order of $M^{0}$ terms

\begin{equation}\label{eq:count_M_0}
\begin{split}
  & \frac{\partial {\alpha}^{(0)}_1}{\partial {t}} + {\mathbf{u}}^{(0)} \cdot \nabla \alpha_1  = 0  \\
  &  \frac{\partial {\alpha}^{(0)}_j  {\rho}^{(0)}_j }{\partial {t}} + \nabla \cdot \left({\alpha}^{(0)}_j {\rho }^{(0)}_j {\mathbf{u}}^{(0)} \right) = 0 \\
  &  \frac{\partial {\rho }^{(0)} {\mathbf{u}}^{(0)} }{\partial {t}} + \nabla \cdot \left( {\rho}^{(0)} {\mathbf{u}}^{(0)} \otimes {\mathbf{u}}^{(0)} \right) + \nabla {p}^{(2)} = 0 \\
  &  \frac{\partial {\rho }^{(0)} {E }^{(0)}}{\partial {t}} + \nabla \cdot \left( \left( {\rho }^{(0)} {E }^{(0)} + {p}^{(0)} \right) {\mathbf{u}}^{(0)} \right) = 0  \\
      &   \frac{\partial  {\alpha}^{(0)}_j  {\rho}^{(0)}_j  {e}^{(0)}_j }{\partial {t}} + \nabla \cdot \left({\alpha}^{(0)}_j  {\rho}^{(0)}_j  {e}^{(0)}_j  {\mathbf{u}}^{(0)} \right) + {\alpha}^{(0)}_j  {p}^{(0)}_j  \nabla \cdot {\mathbf{u}}^{(0)}  = 0 
\end{split}
\end{equation} 
\end{enumerate}

From \eqref{eq:count_M_2} and \eqref{eq:count_M_1} it is evident that leading order $\left( p^{(0)} \right)$ and first order $\left( p^{(1)} \right)$ mean pressure remain uniform in space. Thus mean pressure can be written as

\[p(x, t) = p^{(0)}(t) + p^{(2)} (x, t) M^2 \]

Here, $p^{(0)}(t)$ also contain first order term $\left( p^{(1)} M \right)$. As shown by \citet{pelanti2017low}, the asymptotic results of equilibrium five equation model is similar to the non-equilibrium six-equation model. The only key difference is that in the former system equilibrium pressure $(p = p_1 = p_2)$ will be function of equilibrium Mach number $M_W$. Which is based on equilibrium sound speed $(a_W)$. 

\[p(x, t) = p^{(0)}(t) + p^{(2)} (x, t) M^2_W \]

\section{Low Mach number correction}\label{sec:low_mach}
 As reported by the researchers~\cite{thornber2008improved, rieper2011low, osswald2016l2roe} in the past, scaling the velocity difference with  Mach number can correct the behaviour of pressure in the low Mach number limit. Additionally, it will also reduce the excessive diffusion caused by the Riemann solver. To implement this idea we use reconstructed velocities suggested by the \citet{thornber2008improved}. Expressions for reconstructed normal and tangential velocities can be written as

\begin{equation}\label{eq:reconstruct_velo}
\begin{gathered}
 {u}^{r}_{nL} = \frac{u_{nR} + u_{nL}}{2} - f(M) \frac{u_{nR} - u_{nL}}{2}, \quad {u}^{r}_{tL} = \frac{u_{tR} + u_{tL}}{2} - f(M) \frac{u_{tR} - u_{tL}}{2} \\
 {u}^{r}_{nR} = \frac{u_{nR} + u_{nL}}{2} + f(M) \frac{u_{nR} - u_{nL}}{2}, \quad {u}^{r}_{tR} = \frac{u_{tR} + u_{tL}}{2} + f(M) \frac{u_{tR} - u_{tL}}{2} 
\end{gathered}
\end{equation}
Here, $f(M)$ is function of local Mach number, which can be taken as

\begin{equation}
f(M) = \min\left(1, \max\left(\frac{u^2_{nL} + u^2_{tL}}{a_L}, \frac{u^2_{nR} + u^2_{tR}}{a_R}\right)\right)
\end{equation}

For implementing the low Mach number correction in the Riemann solver, the reconstructed velocities $\left({u}^{r}_{nL}, {u}^{r}_{tL}, {u}^{r}_{nR}, {u}^{r}_{tR} \right)$ will be used in place of left and right velocities $\left( {u}_{nL}, {u}_{tL}, {u}_{nR}, {u}_{tR} \right)$ to compute the fluxes. 
It should be noted that in case of $f(M) = 1.0$ reconstructed left and right velocities in \eqref{eq:reconstruct_velo} become original left and right velocity. 

\section{Asymptotic analysis of semi-discrete equations}\label{sec:asympt_discrete}
To demonstrate the impact of the Low Mach number correction on the pressure scaling of discrete solution, we have conducted asymptotic analysis of the semi-discrete momentum equations. These equations are derived for a  2D quadrilateral mesh using first order scheme with the HLLC solver. In the low Mach number regime, the normal flux vector $\left(\mathbf{F} \left( \mathbf{\hat{Q}}_{L}, \mathbf{\hat{Q}}_{R} \right) \right)$ always lies in the intermediate star region of the HLLC solver, leading to two possible states $\left(\mathbf{F}^{*}_L, \mathbf{F}^{*}_R\right)$ based on the sign of $S^{*}$. Since, the expressions for the fluxes $\mathbf{F}^{*}_L$ and $\mathbf{F}^{*}_R$ are similar, asymptotic expansion for both the cases will be same. Thus, we focused our asymptotic analysis for the case where $\mathbf{F} \left( \mathbf{\hat{Q}}_{L}, \mathbf{\hat{Q}}_{R} \right) = \mathbf{F}^{*}_L$. 

For non-dimensionalisation of semi-discrete equations we have used same reference quantities given in~\eqref{eq:reference}. For representing  expressions involving the difference or sum of variables from cell $i$ and $l$, we used following operators. 

\begin{equation*}
\begin{gathered}
  \Delta_{il}(\phi) = \phi_l - \phi_i \\
   \sigma_{il}(\phi) = \phi_l + \phi_i
\end{gathered}
\end{equation*}  

 With the similar procedure as explained in Section \ref{sec:asympt_contin}, we can get equations for different order of $M$. From these equations we finally arrive on the following results. 

\begin{equation}
 p^{(0)}_i =  p^{(0)} \quad \forall i.
\end{equation} 

\begin{equation}\label{eq:x_mom_-1_1}
 \begin{split}
  & \sum^{4}_{l = 1} \left[ {p^{(1)}_{i} (n_x)_{il}}  + \frac{  {\left({\rho}_i {a}^2_i \right)}^{(0)} \Delta_{il} (p^{(1)}) (n_x)_{il} }{ {a}^{(0)}_i \sigma_{il}( {\rho}^{(0)} {a}^{(0)} )   }  \right] \Delta s_{il} \\
  & - { f(M)} \sum^{4}_{l = 1} \left[ \frac{{\rho}^{(0)}_i {({a}^{(0)}_i)}^2 \left\lbrace \Delta_{il} ({u}^{(0)}_{n}) \Delta_{il}({\rho}^{(0)} {a}^{(0)}) (n_x)_{il} + \sigma_{il}\left({\rho}^{(0)} {a}^{(0)} \right)  \Delta_{il} (u^{(0)}_t ) {(n_y)}_{il} \right\rbrace }{2~ {a}^{(0)}_i \sigma_{il}( {\rho}^{(0)} {a}^{(0)} )}   \right] \Delta s_{il} \\
  & - { f(M)} \sum^{4}_{l = 1} \left[{\rho}^{(0)}_i {a}^{(0)}_i  \frac{ \Delta_{il} (u^{(0)})
  }{2} \right] \Delta s_{il} = 0 
  \end{split}
\end{equation}

\begin{equation}\label{eq:y_mom_-1_1}
 \begin{split}
  & \sum^{4}_{l = 1} \left[{ p^{(1)}_{i} (n_y)_{il}}  + \frac{   {\left({\rho}_i {a}^2_i \right)}^{(0)} \Delta_{il}( p^{(1)}) (n_y)_{il} }{ {a}^{(0)}_i \sigma_{il}( {\rho}^{(0)} {a}^{(0)} )  }  \right] \Delta s_{il} \\
  & - { f(M)} \sum^{4}_{l = 1} \left[ \frac{{\rho}^{(0)}_i {({a}^{(0)}_i)}^2  \left\lbrace \Delta_{il} ({u}^{(0)}_{n}) \Delta_{il}({\rho}^{(0)} {a}^{(0)}) (n_y)_{il} - \sigma_{il}\left({\rho}^{(0)} {a}^{(0)} \right)  \Delta_{il} (u^{(0)}_t ) {(n_x)}_{il} \right\rbrace }{{a}^{(0)}_i \sigma_{il}( {\rho}^{(0)} {a}^{(0)} ) }    \right] \Delta s_{il}  \\
  & - { f(M)} \sum^{4}_{l = 1} \left[{\rho}^{(0)}_i {a}^{(0)}_i   \frac{ \Delta_{il} (v^{(0)})}{2}  \right] \Delta s_{il} = 0
  \end{split}
\end{equation}

For the original HLLC scheme, $f(M)$ is one. After substituting $f(M) = 1$ in above equations \eqref{eq:x_mom_-1_1} and  \eqref{eq:y_mom_-1_1}, it becomes apparent that the first order pressure $(p^{(1)})$ is not constant in space. Consequently, the pressure scaling of the discrete solution obtained by the original HLLC solver differs from that of the continuous system. However, HLLC scheme with the low Mach number correction~\eqref{eq:reconstruct_velo} resolves this issue. In the low Mach number regime $(M<<1)$,  $f(M)$ used in reconstructed velocities~\eqref{eq:reconstruct_velo} becomes the local Mach number. Substituting $f(M) = M$ in \eqref{eq:x_mom_-1_1} and  \eqref{eq:y_mom_-1_1} elevates the coefficient terms to a higher order and as a result we get the following equations.

 \begin{equation}\label{eq:x_mom_-1_2}
 \begin{split}
  & \sum^{4}_{l = 1} \left[ {p^{(1)}_{i} (n_x)_{il}}  + \frac{  {\left({\rho}_i {a}^2_i \right)}^{(0)} \Delta_{il} (p^{(1)}) (n_x)_{il} }{ {a}^{(0)}_i \sigma_{il}( {\rho}^{(0)} {a}^{(0)} )   }  \right] \Delta s_{il}  = 0 
  \end{split}
\end{equation}

\begin{equation}\label{eq:y_mom_-1_2}
 \begin{split}
  & \sum^{4}_{l = 1} \left[{ p^{(1)}_{i} (n_y)_{il}}  + \frac{   {\left({\rho}_i {a}^2_i \right)}^{(0)} \Delta_{il}( p^{(1)}) (n_y)_{il} }{ {a}^{(0)}_i \sigma_{il}( {\rho}^{(0)} {a}^{(0)} )  }  \right] \Delta s_{il} = 0
  \end{split}
\end{equation}

As we can see from the above equations~\eqref{eq:x_mom_-1_2} and~\eqref{eq:y_mom_-1_2}, $\Delta_{il}( p^{(1)})$ should be zero and first order pressure $(p^{(1)})$ should be uniform in space. Thus, with the low Mach number correction HLLC scheme will exhibit the same pressure scaling as the continuous system.

 \section{Results}\label{sec:results}
To demonstrate the effectiveness of the proposed algorithm, we present numerical results for various test cases, including subsonic nozzle flow, dam-break, and low-amplitude sloshing. These cases highlight the algorithm's accuracy in low-Mach-number regimes. Additionally, to assess its robustness and versatility across all speed ranges, we apply the algorithm to a high-speed problem involving shock–helium bubble interaction.

\begin{table}[H]
 \caption{Properties of air and water}
  \centering
  \begin{tabular}{c|c|c|c}
    \hline
      & Density (kg/m3) & $\gamma$ & $\pi$  \\  \hline
    Air  & 1 &   1.4 & 0 \\
    Water & 1000 &  4.4 & $6 \times 10^8$ \\
     \hline
  \end{tabular}
  \label{tab:properties}
\end{table}

\subsection{Subsonic flow in symmetric nozzle}

\begin{figure}[H]
           \centering
             \includegraphics[width =0.8\textwidth]{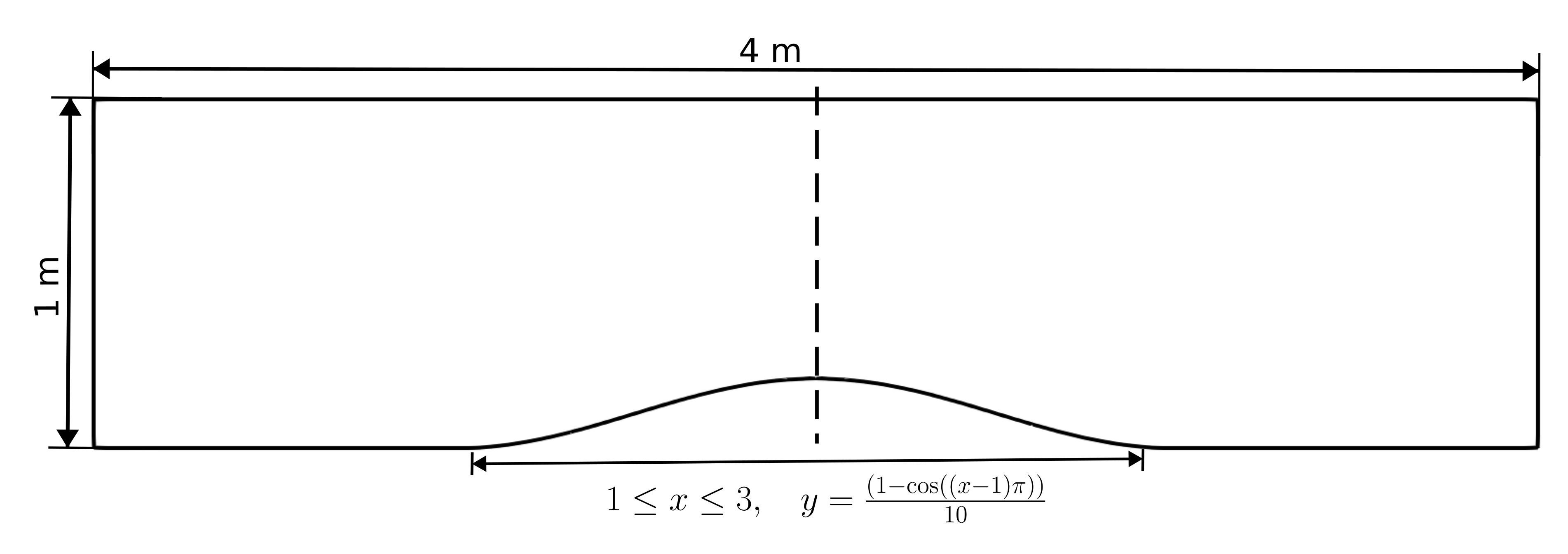}
           \caption{Nozzle geometry}
           \label{fig:nozzle_geo}
\end{figure}

\begin{figure}[H]
           \centering
             \includegraphics[width =1\textwidth]{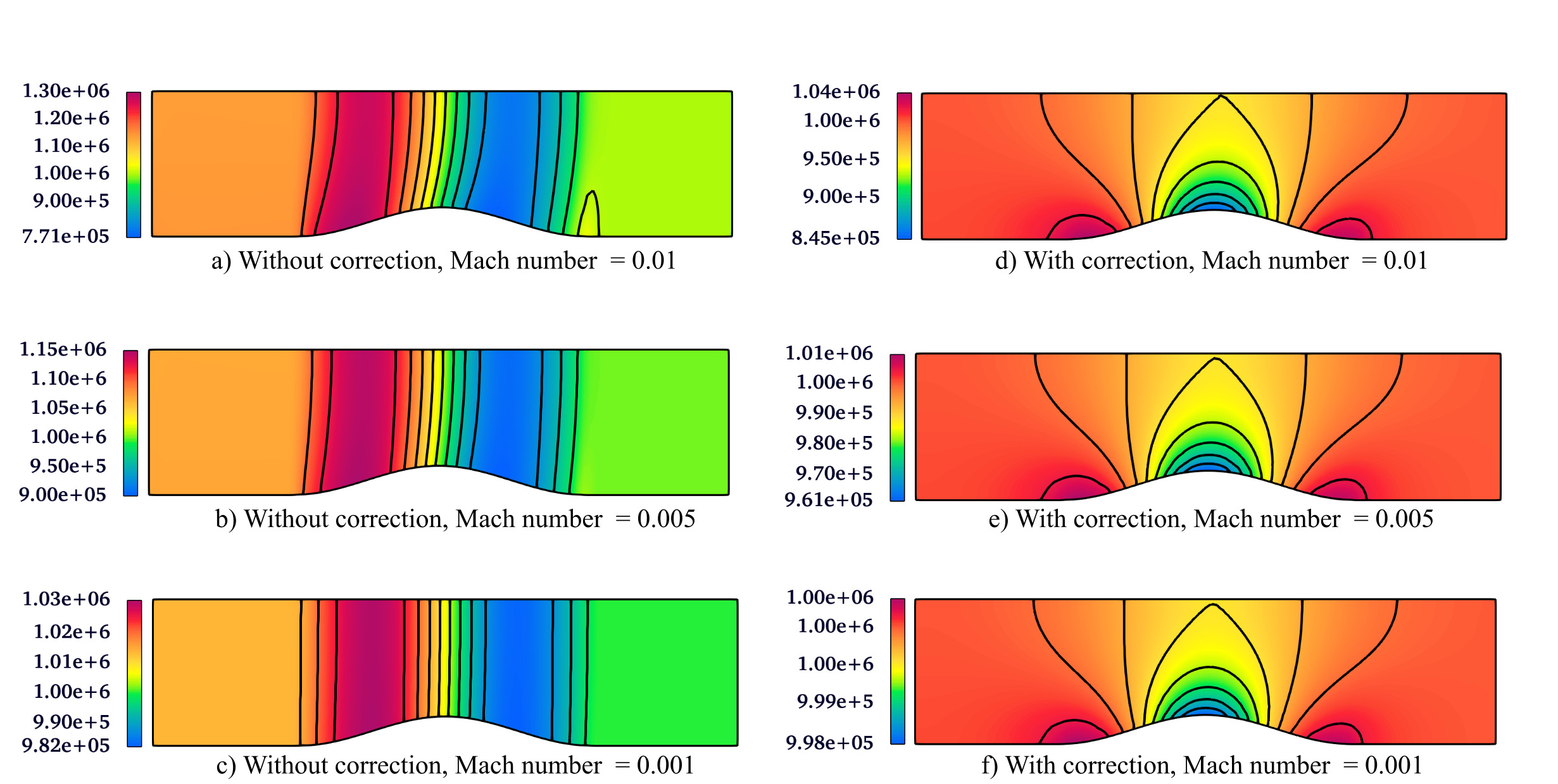}
           \caption{ Pressure contours of subsonic nozzle problem for different Mach numbers.}
           \label{fig:Contours}
\end{figure}

\begin{figure}[H]
           \centering
             \includegraphics[width = 1\textwidth]{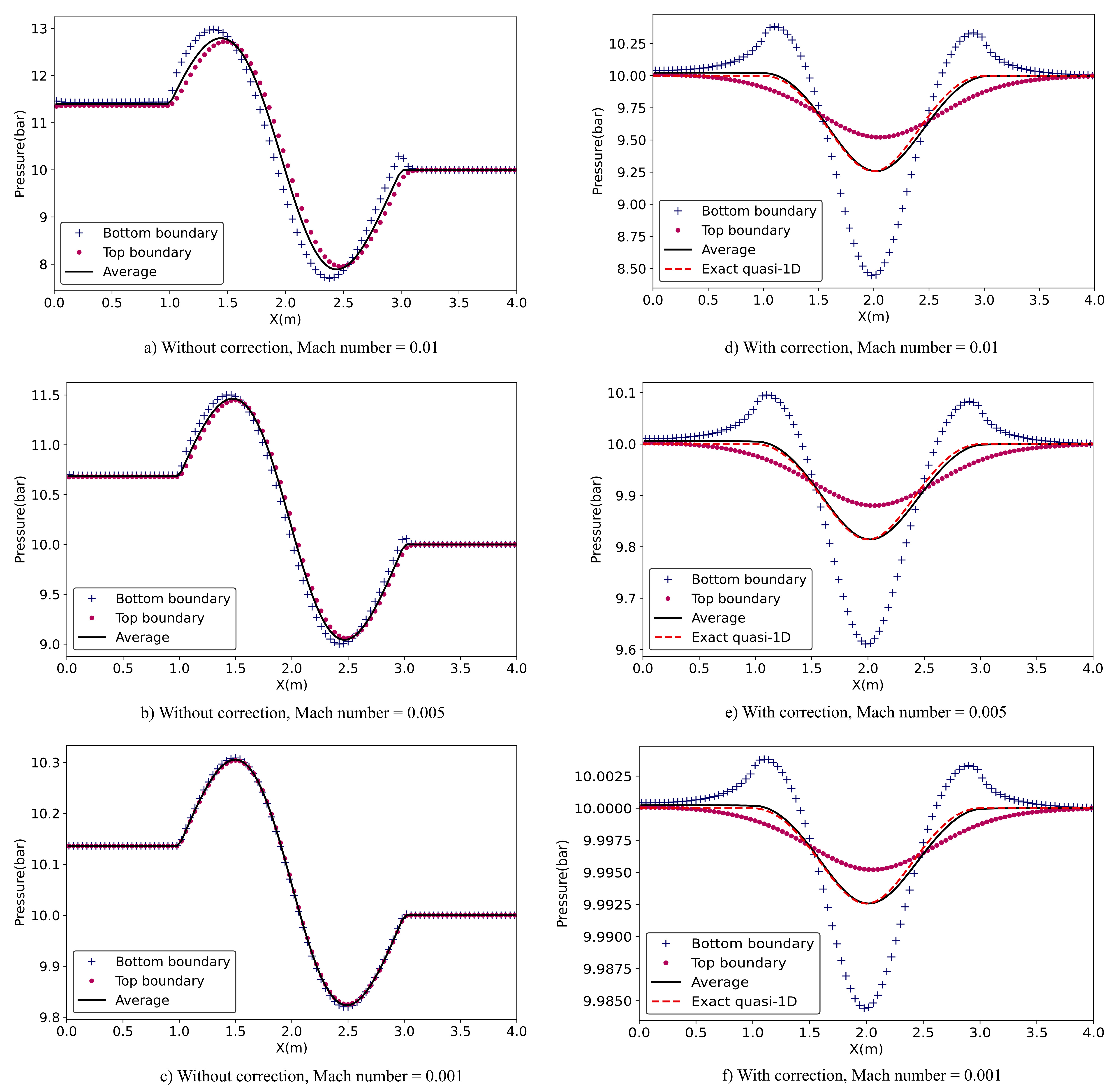}
           \caption{Pressure plot of subsonic nozzle problem for different inlet Mach number.}
           \label{fig:Plots}
\end{figure}

\begin{figure}[H]
           \centering
             \includegraphics[scale =  0.6]{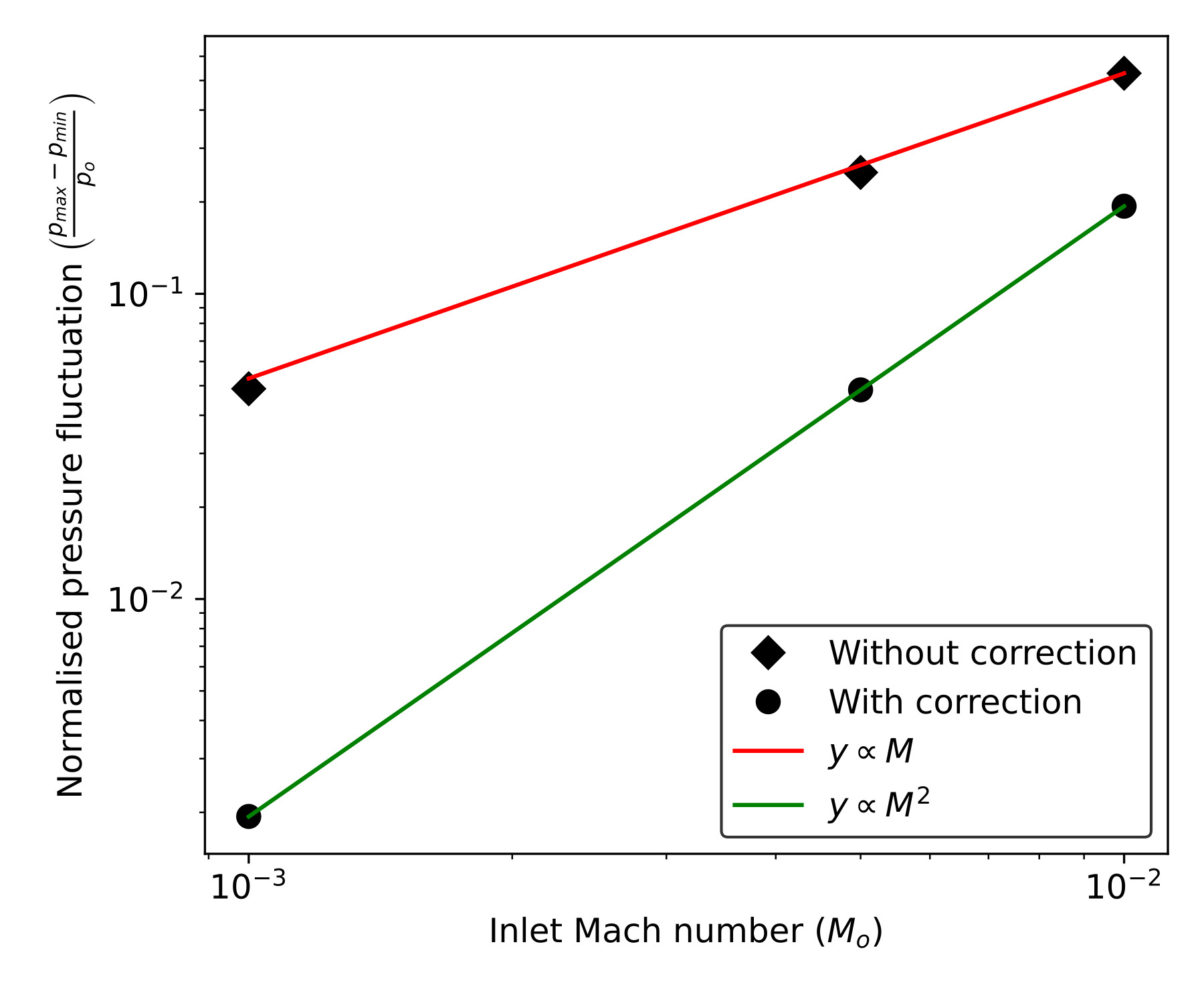}
           \caption{ Log-log plot of computed normalised pressure fluctuation vs inlet Mach number.}
           \label{fig:fluctuation_Plot}
\end{figure}

We considered a two-phase nozzle problem, in which a water-air mixture flows inside a two-dimensional nozzle at subsonic Mach number. The nozzle geometry used in the problem, shown in Figure 2, is taken from a similar test case presented in the literature~\cite{murrone2008, pelanti2017low}. The computational domain is discretized with $100\times25$ quadrilateral cells. To avoid interference from limiters and spatial reconstruction, only first-order results are investigated.  In this test case, a water-air mixture with $\alpha_{air} = 10^{-3}$ is considered at inlet. Water and air properties are provided in \tbl{tab:properties}. A fixed pressure, $p_o = 10^6$ Pa is imposed on the outlet. The inflow velocity is based on the selected Mach number. In this work, three sets of numerical experiments corresponds to Mach number $M_0 = 0.01, 0.005, 0.001$ are performed. The results obtained with and without correction are shown \fig{fig:Contours} and \fig{fig:Plots}. As the nozzle is symmetric about the middle vertical axis, the flow inside the nozzle should also be symmetric. However, if we observe pressure contours from Figure \ref{fig:Contours} a) to
\ref{fig:Contours} c), the results obtained by standard HLLC scheme is not symmetric. This unphysical behaviour is rectified after applying the proposed correction, as it can be observed from \ref{fig:Contours} d) to
\ref{fig:Contours} f). 
In addition to the qualitative comparison, quantitative results are also presented in the Figure \ref{fig:Plots}. The plots shown in the Figure \ref{fig:Plots} contain pressure profile of bottom and top wall of the nozzle, as well as the
average pressure over the height. Additionally, average pressure curve also compared with the exact solution of quasi one-dimensional flow for the results obtained with the correction. The comparison presented in \fig{fig:Plots} indicates that the proposed correction effectively resolves the unphysical behaviour observed in the existing numerical scheme, and the results align well with the exact solution.
From the asymptotic analysis of continuous model, we know that pressure fluctuation in low Mach number limit is order of square of Mach number. To confirm this behaviour in the numerical results, the normalised pressure fluctuation is plotted against inlet Mach number in the Figure~\ref{fig:fluctuation_Plot}. From the log-log plot (Figure~\ref{fig:fluctuation_Plot}), it becomes apparent that pressure scaling in the numerical results obtained using the proposed correction follows the correct behaviour. 
\subsection{Dam break problem}
\begin{figure}[H]
           \centering
             \includegraphics[ width =0.6\textwidth]{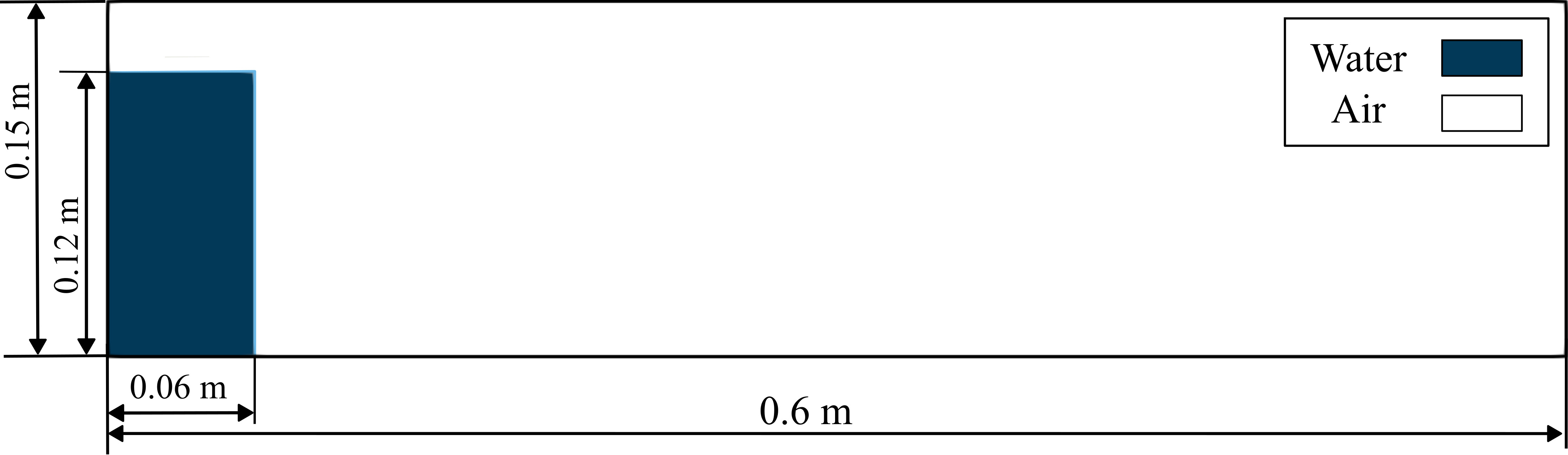}
           \caption{Initial setup for the dam break problem.}
           \label{fig:dam_break_initial}
\end{figure}

\begin{figure}[H]
           \centering
             \includegraphics[ width =1\textwidth]{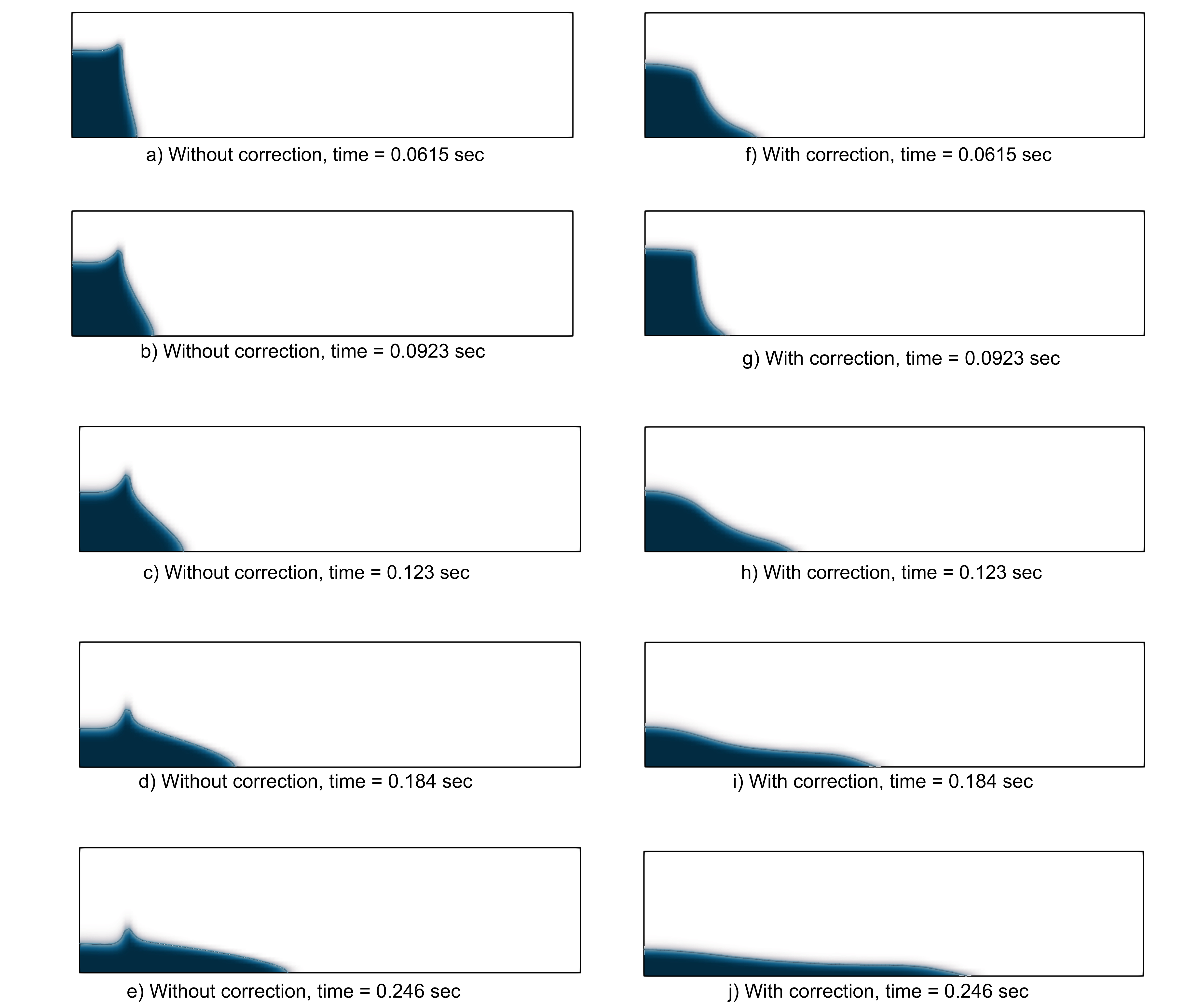}
           \caption{Volume fraction contours of the dam break problem at different time.}
           \label{fig:comparison_dam_break}
\end{figure}

\begin{figure}[H]
           \centering
             \includegraphics[ width =1\textwidth]{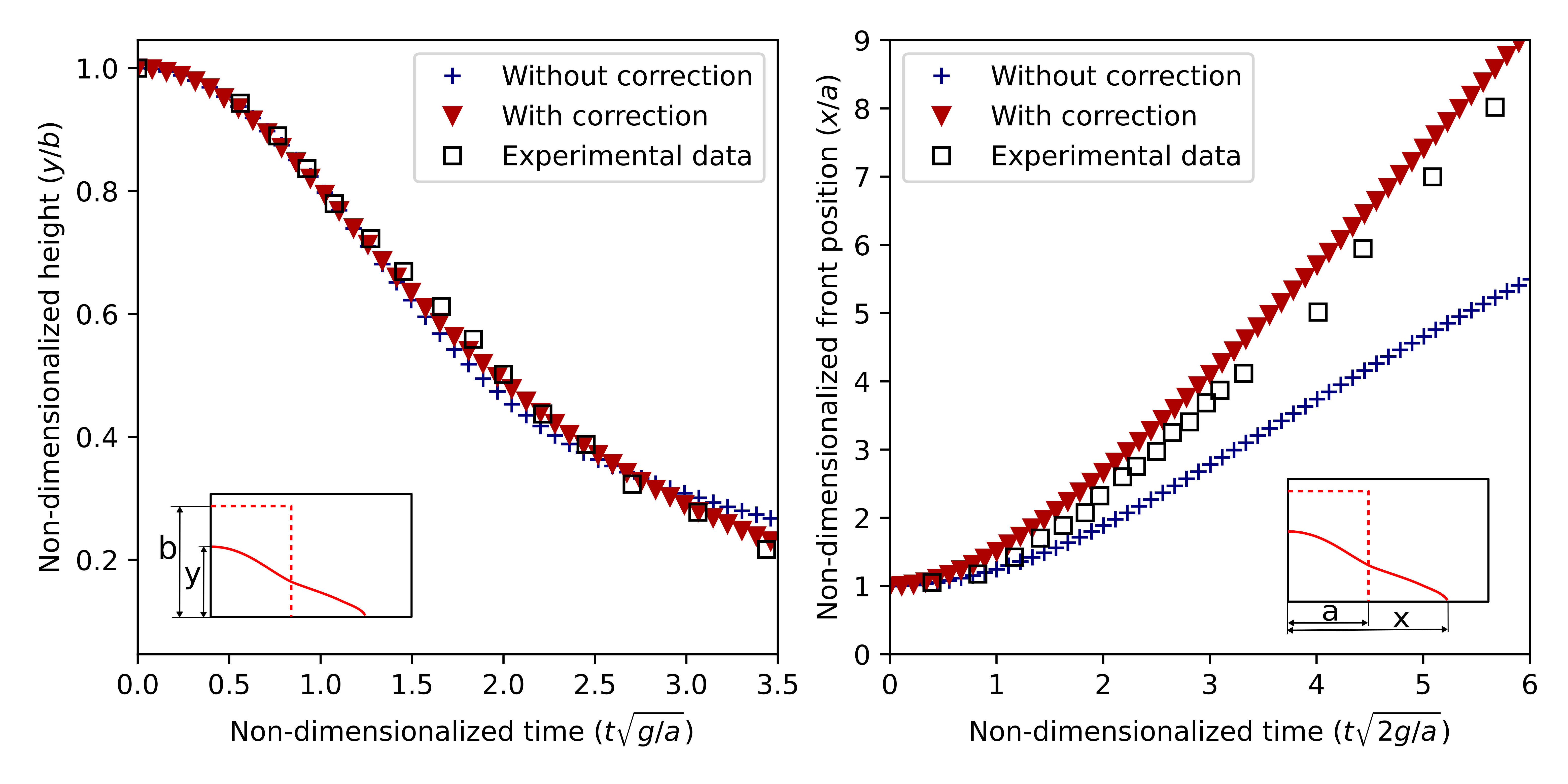}
           \caption{\centering Comparison between numerical solutions obtained with and without low Mach
number correction and experimental data for dam break problem.}
           \label{fig:comparison_plot_dam_break}
\end{figure}

This is a standard test case used to test numerical methods for resolved interface flows at low Mach number \citep{murrone2005five, murrone2008}. The problem involves a water column inside a closed domain filled with air. The gravitational force causes the water column to collapse. The initial setup for the problem is presented in \fig{fig:dam_break_initial}. The numerical results are obtained using second-order scheme on $120 \times 30$ structured mesh.
\fig{fig:comparison_dam_break} presents a comparison of the results obtained with and without the proposed correction to the numerical scheme, illustrating the improvements achieved after applying the correction. From the \fig{fig:comparison_dam_break} a) to \ref{fig:comparison_dam_break} e), it can be observed that the results obtained using standard scheme are unphysical. Specifically, the water-air interface is not smooth and exhibits an unphysical peak at the top right corner. This peak remains attached to the interface even after the water column becomes almost flat. However, from \fig{fig:comparison_dam_break} f) to \ref{fig:comparison_dam_break} j), it is evident that the unphysical behaviour in the numerical results vanishes after applying the low Mach correction, resulting in a smooth water-air interface. \fig{fig:comparison_plot_dam_break} shows quantitative comparison of the numerical results against experimental data \cite{martin1952experimental}. The comparison is made by plotting the non-dimensionalized height $(y/b)$ and non-dimensionalized front position $(x/a)$. In this study, initial width $(a)$ and height $(b)$ of the water column are taken as, 0.06 m and 0.12 m respectively. From the  \fig{fig:comparison_plot_dam_break}, it is apparent that difference between experimental data and numerical results is reduced after applying low Mach number correction.
\subsection{Low amplitude sloshing}
This is a classical test problem \cite{yang2014upwind, bhat2019contact, parameswaran2023stable}, where a water-air interface oscillates under the influence of the gravity. Since, the experiment is conducted without considering any viscous and surface tension effects, we expect the oscillations to be continued without any damping.   
\begin{figure}[H]
           \centering
             \includegraphics[width =0.35\textwidth]{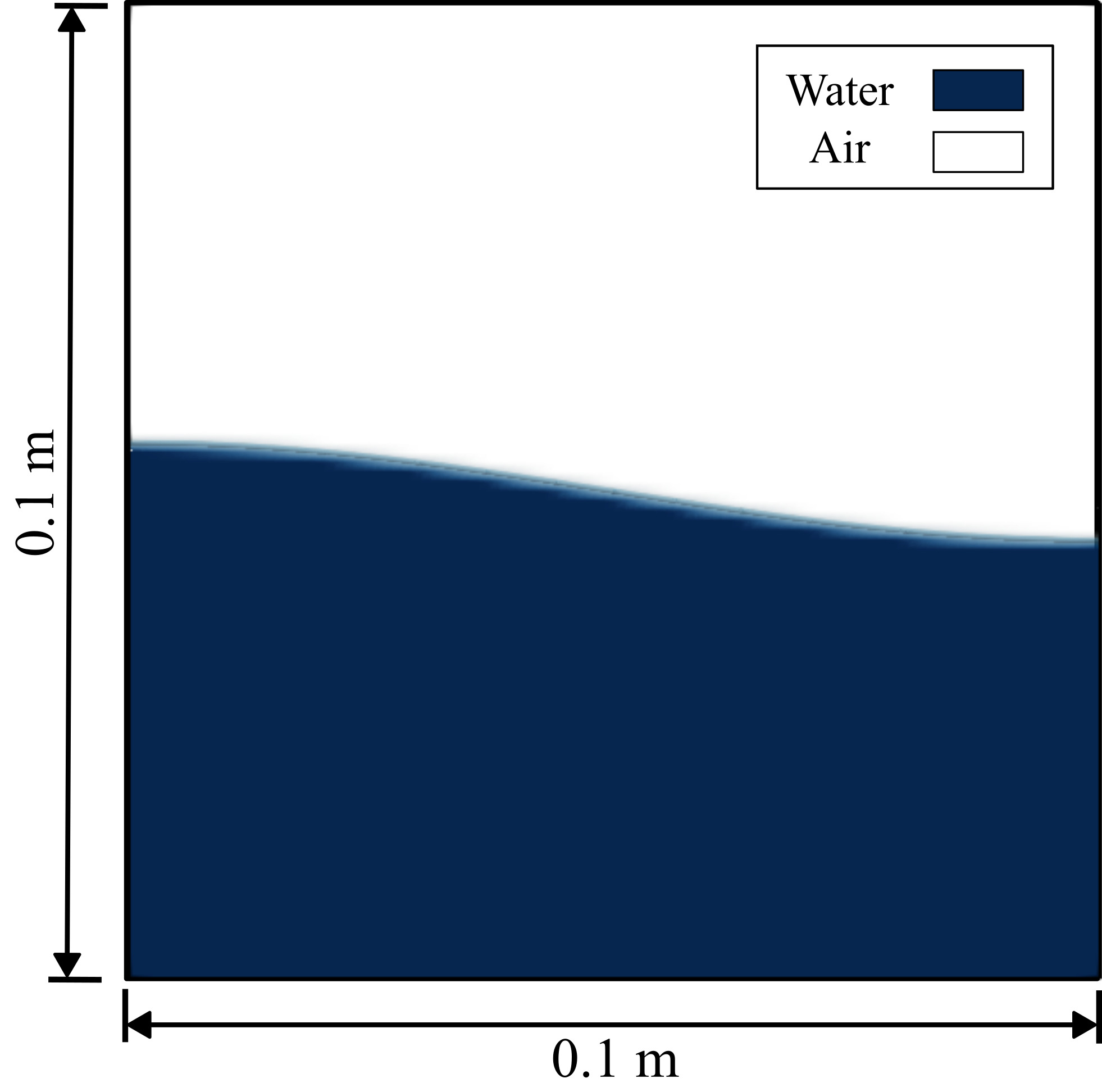}
           \caption{Initial setup for the Low amplitude sloshing problem.}
           \label{fig:initial_sloshing}
\end{figure}
The initial configuration for the test case is illustrated in \fig{fig:initial_sloshing}, where water and air are contained separately in the square shaped box with side length, $L = 0.1$ m. The two-dimensional computational space is discretised with $100 \times 100$ uniform cells and numerical results are obtained using second order scheme. Initially, water-air interface follows half of cosine curve, given by the equation, $y(x) = 0.05 + 0.005 \cos(\pi x/L)$. Under the influence of gravity, interface begins to oscillate. The time period for the first mode of oscillation is determined using the following expression \cite{tadjbakhsh1960standing},
$$ P =  2 \pi \sqrt{g k \tanh(kh)}. $$
Where, $h$ denotes average depth of water and $k = \frac{2 \pi}{\lambda}$ is the wave number. For this problem, time period of the first mode calculated to be 0.3739 sec. 

\begin{figure}[H]
           \centering
           \includegraphics[ width = 0.7\textwidth]{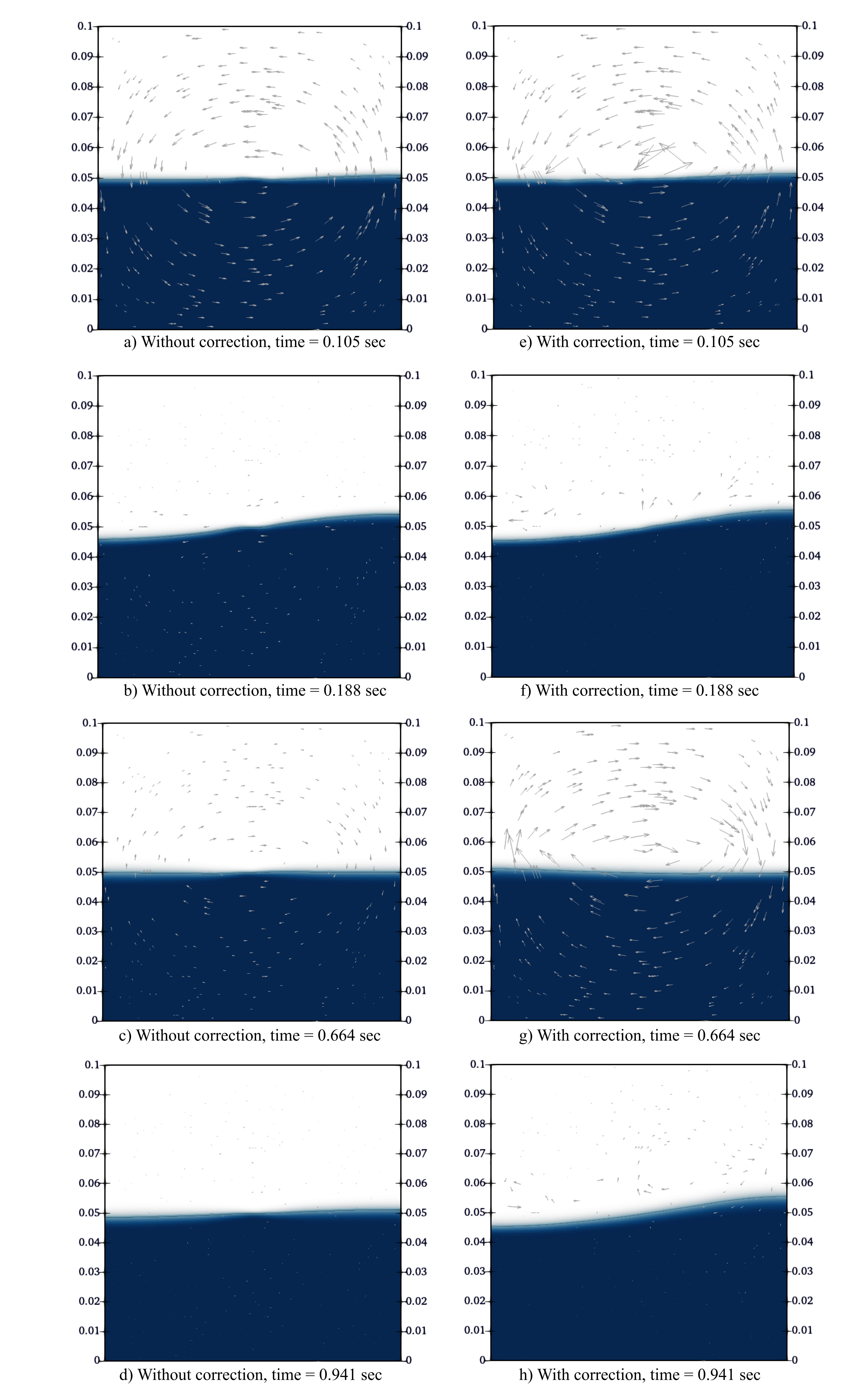}
           \caption{Volume fraction contours and velocity field at different time intervals for sloshing problem.}
           \label{fig:comparison_sloshing}
\end{figure}

\begin{figure}[H]
           \centering
             \includegraphics[width = 1\textwidth]{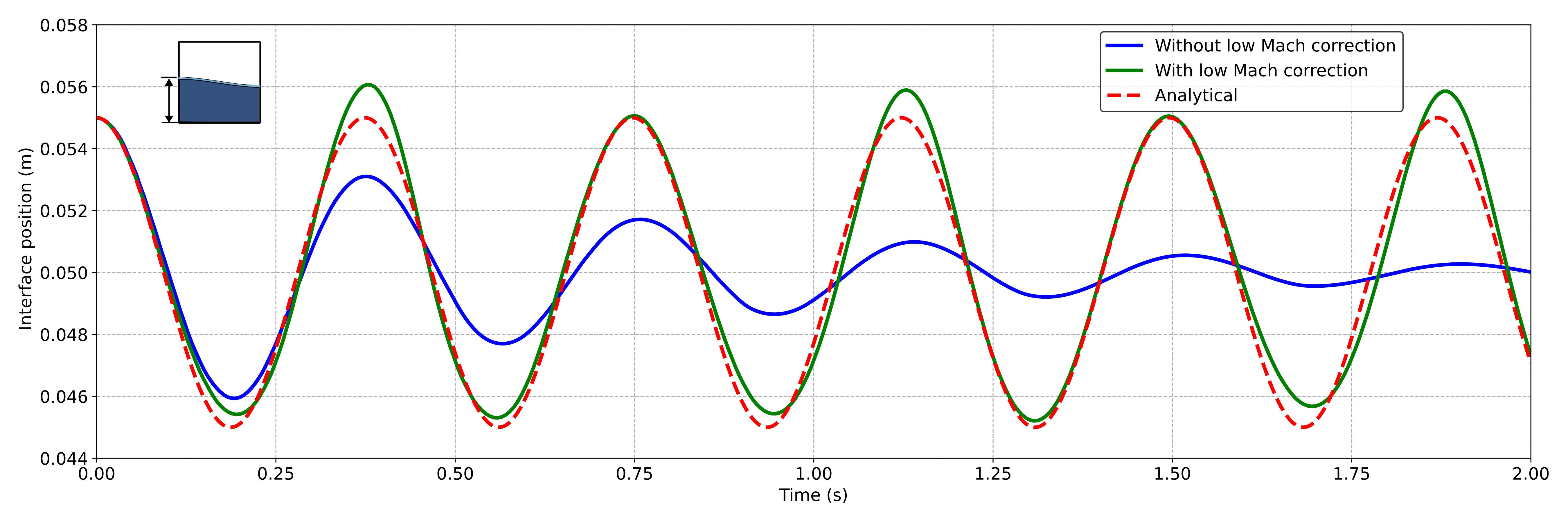}
           \caption{\centering Comparison of low amplitude sloshing results obtained with and without low Mach number correction and analytical solution for first mode of oscillation.}
           \label{fig:sloshing_plot}
\end{figure}

A series of images of volume fraction contour with the velocity field are plotted in \fig{fig:comparison_sloshing}. The numerical results obtained by the standard HLLC scheme are presented on the left side of the \fig{fig:comparison_sloshing}, while right side shows results obtained with the correction. From the sequence of images in  \fig{fig:comparison_sloshing}, one can observe that as the time progresses, the amplitude of oscillation decreases in the results obtained by the standard HLLC scheme. In contrast, in the results obtained with the proposed correction, the interface location at the extrema of oscillation is almost same at time = 0.188 sec (\fig{fig:comparison_sloshing} f)) and time = 0.941 sec (\fig{fig:comparison_sloshing} h)). This implies that numerical viscosity is reduced when the proposed correction is applied to the scheme. To assess the numerical results, the water-air interface location at the left wall of the container is plotted in \fig{fig:sloshing_plot}. The plot presented in the \fig{fig:sloshing_plot} shows that interface motion slows down over time and nearly ceases at the end for the result obtained by the standard HLLC scheme. However, if we use the proposed correction, the oscillation is maintained even at the end of the plot, and the peak of the first mode aligns with the analytical solution.
\subsection{Shock-bubble interaction}
\begin{figure}[H]
           \centering
             \includegraphics[width = 0.8\textwidth]{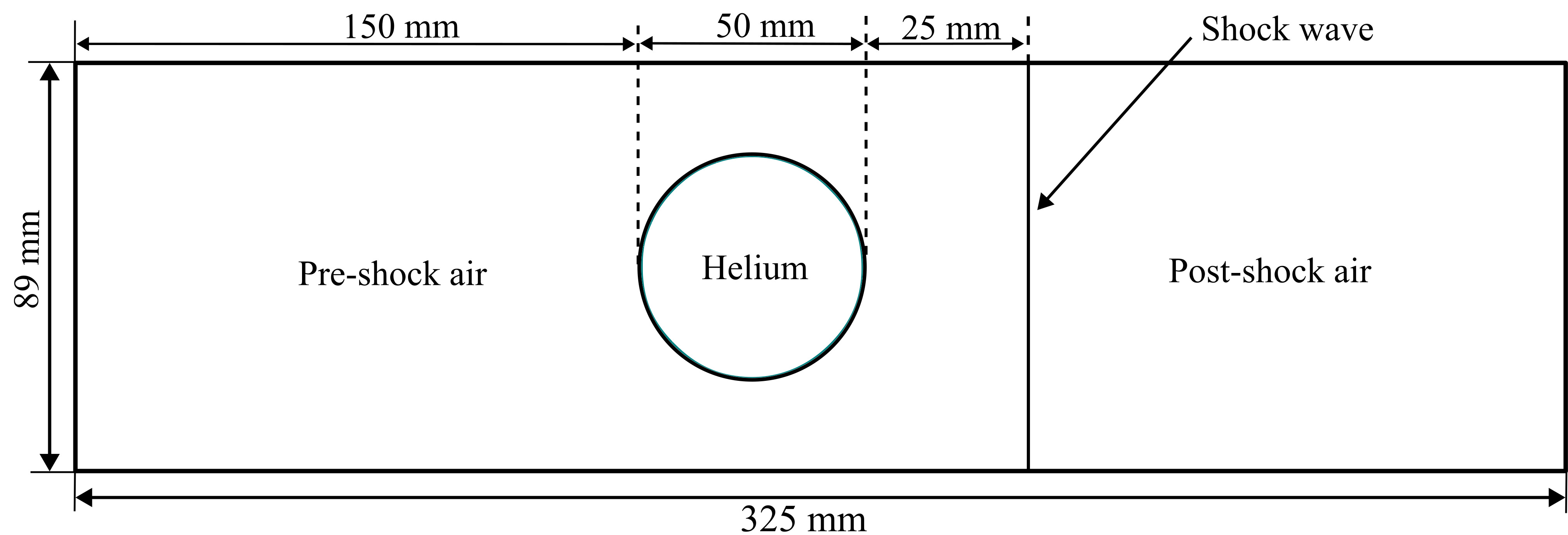}
           \caption{\centering Computational domain for shock-bubble interaction problem.}
\label{fig:Initial_shock_bubble}
\end{figure}

In order to test the numerical scheme with proposed correction for high-speed flows and shock waves, a standard problem of shock-helium bubble interaction \cite{quirk1996dynamics, terashima2009front, yoo2018numerical, nguyen2022fully} is solved.
\begin{table}[H]
 \caption{Initial parameters for shock-bubble problem}
  \centering
  \begin{tabular}{c|c|c|c|c|c}
    \hline
      & $\gamma$ & $\pi$  & Density $(kg/m^3)$ & Velocity $(m/s)$ & Pressure $(Pa)$ \\  
     \hline
    Pre-shock air  &  1.4 & 0 & 1.4 & (0, 0) & 100000\\
    Post-shock air  &    1.4 & 0 & 1.92691 & (-114.42, 0) & 156980\\
    Helium bubble &   1.648 & 0 & 0.25463 & (0, 0) & 100000 \\
     \hline 
  \end{tabular}
  \label{tab:properties_shock_bubble}
\end{table}
At initial stage of the simulation, a stationary helium bubble is placed before a moving shock wave. The Initial configuration, including geometrical details of the problem is shown in \fig{fig:Initial_shock_bubble}. The parameters for different regions, as shown in \fig{fig:Initial_shock_bubble} are provided in \tbl{tab:properties_shock_bubble}. A uniform mesh of $650 \times 180$ cells \cite{nguyen2022fully} are used for the discretization. All the primitive variables, except for volume fraction $(\alpha_1)$, are reconstructed using the SDWLS method. For the volume fraction $(\alpha_1)$, the overbee limiter \cite{overbee} is used.

\begin{figure}[H]
           \centering
             \includegraphics[scale = 0.4]{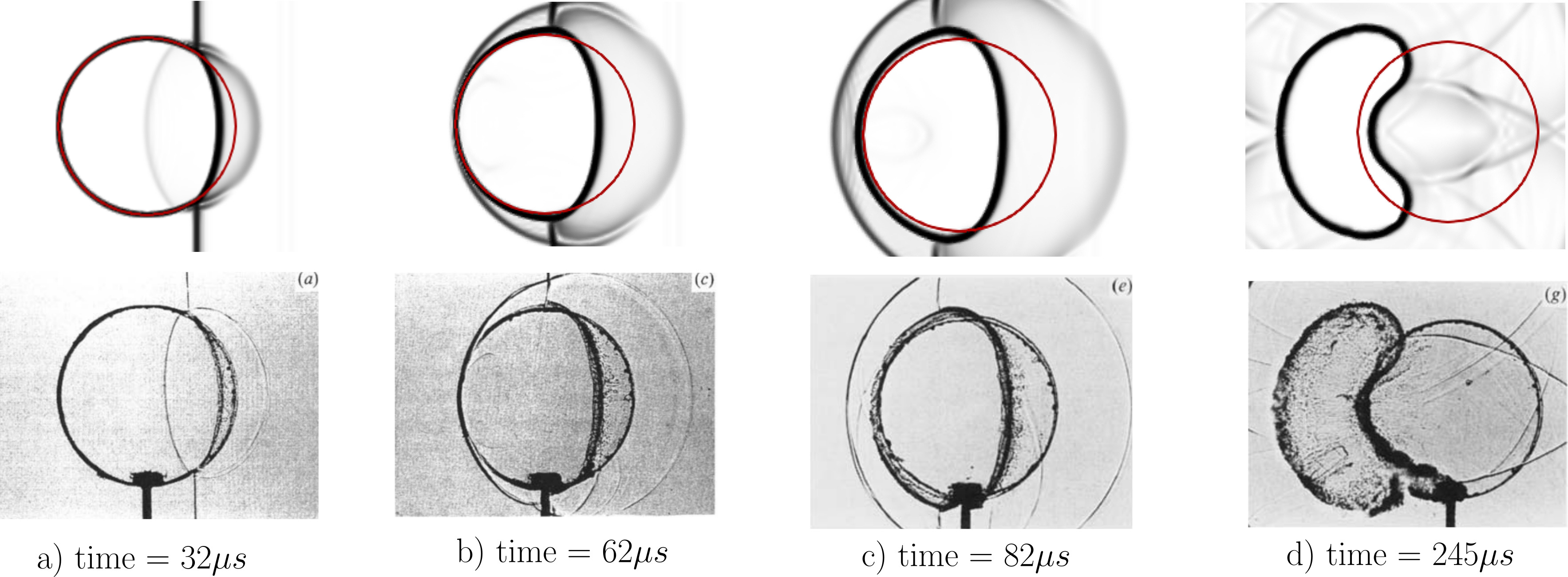}
           \caption{\centering Comparison of numerical results (top row) with experimental images \cite{haas1987interaction} (bottom row).}
\label{fig:shock_bubble_images}
\end{figure}

A set of schlieren images obtained from the numerical results at various time intervals are presented along with the experimental images \cite{haas1987interaction} in \fig{fig:shock_bubble_images}. The collision of shock wave with stationary bubble, deforms and accelerates the bubble, leading to the reflection and transmission of shock wave. A quantitative comparison in \fig{fig:shock_bubble_images} demonstrate that the numerical scheme with the proposed correction accurately captures all the features observed in the experimental images.  

\begin{figure}[H]
           \centering
             \includegraphics[scale = 0.6]{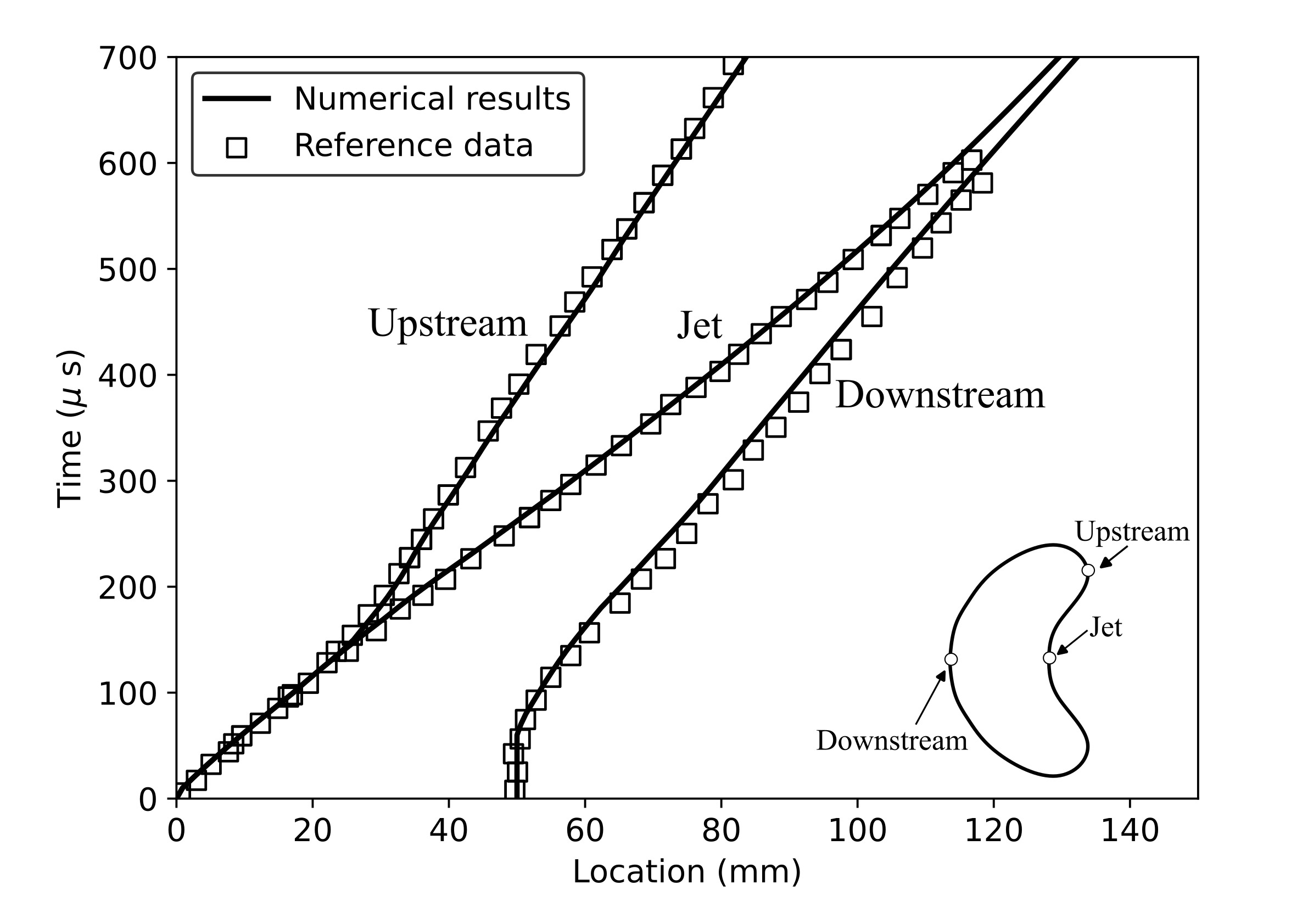}
           \caption{\centering Trajectory of three characteristic points on space-time graph.}
\label{fig:shock_bubble_plot}
\end{figure}

For further validation, the trajectories of three typical points, jet, upstream and downstream on the bubble are tracked and plotted in \fig{fig:shock_bubble_plot}. The space-time curves of these points are also compared with the  front tracking method result \cite{terashima2009front}. In both \fig{fig:shock_bubble_images} and \fig{fig:shock_bubble_plot}, starting time is considered from the moment when the shock wave touches the helium bubble. The comparison shows that numerical results obtained with the proposed correction are in good agreement with the standard reference \cite{terashima2009front} solution.


\section{Conclusion}\label{sec:conclusion}
We presented a numerical scheme capable of solving multiphase problems across all Mach numbers. The six-equation model with instantaneous relaxation was employed for multiphase flow computations. Conventional numerical approaches for solving this model often produce inaccurate results at low Mach numbers. To address this issue, we incorporated the velocity reconstruction formula proposed by Thornber et al. into the HLLC Riemann solver. Unlike preconditioning methods, the proposed correction does not impose a restrictive time step size in explicit schemes. Additionally, it avoids the global cutoff Mach number issue by utilizing local Mach number scaling for velocity differences.

The proposed correction was validated through asymptotic analysis on both the continuous model and its discretized form. The asymptotic expansion results confirm that the pressure scaling of the discrete solution aligns with the expected behavior after applying the correction. We demonstrated its effectiveness through various low-Mach-number multiphase problems, including subsonic nozzle flow, dam-break, and low-amplitude sloshing. A series of nozzle experiments at different Mach numbers further validated that the corrected numerical results adhere to the correct pressure scaling. Comparisons with analytical solutions and experimental data highlight the effectiveness of the correction, as the numerical results closely match reference solutions.

\bibliographystyle{elsarticle-num-names} 

\bibliography{ref}

\end{document}